\documentstyle[11pt]{article}

\newcommand{\Label}{\label}

\begin{document}

\baselineskip16pt

\newtheorem{definition}{Definition $\!\!$}
\newtheorem{prop}[definition]{Proposition $\!\!$}
\newtheorem{lem}[definition]{Lemma $\!\!$}
\newtheorem{corollary}[definition]{Corollary $\!\!$}
\newtheorem{theorem}[definition]{Theorem $\!\!$}
\newtheorem{example}[definition]{\it Example $\!\!$}
\newtheorem{remark}[definition]{Remark $\!\!$}

\newcommand{\nc}[2]{\newcommand{#1}{#2}}
\newcommand{\rnc}[2]{\renewcommand{#1}{#2}}

\nc{\bpr}{\begin{prop}}
\nc{\bth}{\begin{theorem}}
\nc{\ble}{\begin{lem}}
\nc{\bco}{\begin{corollary}}
\nc{\bre}{\begin{remark}}
\nc{\bex}{\begin{example}}
\nc{\bde}{\begin{definition}}
\nc{\ede}{\end{definition}}
\nc{\epr}{\end{prop}}
\nc{\ethe}{\end{theorem}}
\nc{\ele}{\end{lem}}
\nc{\eco}{\end{corollary}}
\nc{\ere}{\hfill\mbox{$\Diamond$}\end{remark}}
\nc{\eex}{\end{example}}
\nc{\epf}{\hfill\mbox{$\Box$}}
\nc{\ot}{\otimes}
\nc{\bsb}{\begin{Sb}}
\nc{\esb}{\end{Sb}}
\nc{\ct}{\mbox{${\cal T}$}}
\nc{\ctb}{\mbox{${\cal T}\sb B$}}
\nc{\bcd}{\[\begin{CD}}
\nc{\ecd}{\end{CD}\]}
\nc{\ba}{\begin{array}}
\nc{\ea}{\end{array}}
\nc{\bea}{\begin{eqnarray}}
\nc{\eea}{\end{eqnarray}}
\nc{\be}{\begin{enumerate}}
\nc{\ee}{\end{enumerate}}
\nc{\beq}{\begin{equation}}
\nc{\eeq}{\end{equation}}
\nc{\bi}{\begin{itemize}}
\nc{\ei}{\end{itemize}}
\nc{\kr}{\mbox{Ker}}
\nc{\te}{\!\ot\!}
\nc{\pf}{\mbox{$P\!\sb F$}}
\nc{\pn}{\mbox{$P\!\sb\nu$}}
\nc{\bmlp}{\mbox{\boldmath$\left(\right.$}}
\nc{\bmrp}{\mbox{\boldmath$\left.\right)$}}
\rnc{\phi}{\mbox{$\varphi$}}
\nc{\LAblp}{\mbox{\LARGE\boldmath$($}}
\nc{\LAbrp}{\mbox{\LARGE\boldmath$)$}}
\nc{\Lblp}{\mbox{\Large\boldmath$($}}
\nc{\Lbrp}{\mbox{\Large\boldmath$)$}}
\nc{\lblp}{\mbox{\large\boldmath$($}}
\nc{\lbrp}{\mbox{\large\boldmath$)$}}
\nc{\blp}{\mbox{\boldmath$($}}
\nc{\brp}{\mbox{\boldmath$)$}}
\nc{\LAlp}{\mbox{\LARGE $($}}
\nc{\LArp}{\mbox{\LARGE $)$}}
\nc{\Llp}{\mbox{\Large $($}}
\nc{\Lrp}{\mbox{\Large $)$}}
\nc{\llp}{\mbox{\large $($}}
\nc{\lrp}{\mbox{\large $)$}}
\nc{\lbc}{\mbox{\Large\boldmath$,$}}
\nc{\lc}{\mbox{\Large$,$}}
\nc{\Lall}{\mbox{\Large$\forall$}}
\nc{\bc}{\mbox{\boldmath$,$}}
\rnc{\epsilon}{\varepsilon}
\rnc{\ker}{\mbox{\em Ker}}
\nc{\ra}{\rightarrow}
\nc{\ci}{\circ}
\nc{\cc}{\!\ci\!}
\nc{\T}{\mbox{\sf T}}
\nc{\can}{\mbox{\em\sf T}\!\sb R}
\nc{\cnl}{$\mbox{\sf T}\!\sb R$}
\nc{\lra}{\longrightarrow}
\nc{\M}{\mbox{Map}}
\nc{\imp}{\Rightarrow}
\rnc{\iff}{\Leftrightarrow}
\nc{\bmq}{\cite{bmq}}
\nc{\ob}{\mbox{$\Omega\sp{1}\! (\! B)$}}
\nc{\op}{\mbox{$\Omega\sp{1}\! (\! P)$}}
\nc{\oa}{\mbox{$\Omega\sp{1}\! (\! A)$}}
\nc{\inc}{\mbox{$\,\subseteq\;$}}
\nc{\de}{\mbox{$\Delta$}}
\nc{\spp}{\mbox{${\cal S}{\cal P}(P)$}}
\nc{\dr}{\mbox{$\Delta_{R}$}}
\nc{\dsr}{\mbox{$\Delta_{\cal R}$}}
\nc{\m}{\mbox{m}}
\nc{\0}{\sb{(0)}}
\nc{\1}{\sb{(1)}}
\nc{\2}{\sb{(2)}}
\nc{\3}{\sb{(3)}}
\nc{\4}{\sb{(4)}}
\nc{\5}{\sb{(5)}}
\nc{\6}{\sb{(6)}}
\nc{\7}{\sb{(7)}}
\nc{\hsp}{\hspace*}
\nc{\nin}{\mbox{$n\in\{ 0\}\!\cup\!{\Bbb N}$}}
\nc{\al}{\mbox{$\alpha$}}
\nc{\bet}{\mbox{$\beta$}}
\nc{\ha}{\mbox{$\alpha$}}
\nc{\hb}{\mbox{$\beta$}}
\nc{\hg}{\mbox{$\gamma$}}
\nc{\hd}{\mbox{$\delta$}}
\nc{\he}{\mbox{$\varepsilon$}}
\nc{\hz}{\mbox{$\zeta$}}
\nc{\hs}{\mbox{$\sigma$}}
\nc{\hk}{\mbox{$\kappa$}}
\nc{\hm}{\mbox{$\mu$}}
\nc{\hn}{\mbox{$\nu$}}
\nc{\la}{\mbox{$\lambda$}}
\nc{\hl}{\mbox{$\lambda$}}
\nc{\hG}{\mbox{$\Gamma$}}
\nc{\hD}{\mbox{$\Delta$}}
\nc{\th}{\mbox{$\theta$}}
\nc{\Th}{\mbox{$\Theta$}}
\nc{\ho}{\mbox{$\omega$}}
\nc{\hO}{\mbox{$\Omega$}}
\nc{\hp}{\mbox{$\pi$}}
\nc{\hP}{\mbox{$\Pi$}}
\nc{\bpf}{{\it Proof.~~}}
\nc{\slq}{\mbox{$A(SL\sb q(2))$}}
\nc{\fr}{\mbox{$Fr\llp A(SL(2,\IC))\lrp$}}
\nc{\slc}{\mbox{$A(SL(2,\IC))$}}
\nc{\af}{\mbox{$A(F)$}}
\rnc{\widetilde}{\tilde}
\nc{\qdt}{quantum double torus}
\nc{\aqdt}{\mbox{$A(DT^2_q)$}}
\nc{\dtq}{\mbox{$DT^2_q$}}
\nc{\uc}{\mbox{$U(2)$}}
\nc{\uq}{\mbox{$U_{q^{-1},q}(2)$}}
\nc{\qg}{quantum group}
\nc{\qpb}{quantum principal bundle}
\nc{\hop}{Hopf algebra}

\def\esl{{\mbox{$E\sb{\frak s\frak l (2,{\Bbb C})}$}}}
\def\esu{{\mbox{$E\sb{\frak s\frak u(2)}$}}}
\def\zf{{\mbox{${\Bbb Z}\sb 4$}}}
\def\zt{{\mbox{$2{\Bbb Z}\sb 2$}}}
\def\ox{{\mbox{$\Omega\sp 1\sb{\frak M}X$}}}
\def\oxh{{\mbox{$\Omega\sp 1\sb{\frak M-hor}X$}}}
\def\oxs{{\mbox{$\Omega\sp 1\sb{\frak M-shor}X$}}}
\def\Fr{\mbox{Fr}}
\def\gal{-Galois extension}
\def\hge{Hopf-Galois extension}
\def\cge{coalgebra-Galois extension}
\def\pge{$\psi$-Galois extension}
\def\ta{\tilde a}
\def\tb{\tilde b}
\def\tc{\tilde c}
\def\td{\tilde d}
\def\st{\stackrel}

\def\inbar{\,\vrule height1.5ex width.4pt depth0pt}
\def\IC{{\Bbb C}}
\def\IZ{{\Bbb Z}}
\def\IN{{\Bbb N}}
\def\otc{\otimes_{\IC}}
\def\ra{\rightarrow}
\def\ota{\otimes_ A}
\def\otza{\otimes_{ Z(A)}}
\def\otc{\otimes_{\IC}}
\def\h{\rho}
\def\x{\zeta}
\def\th{\theta}
\def\s{\sigma}
\def\t{\tau}

\title{
\vspace*{-15mm}
{
\large\bf 
A NOTE ON FIRST ORDER DIFFERENTIAL
CALCULUS ON QUANTUM PRINCIPAL BUNDLES\\
}
\author{
{\sc Piotr M.~Hajac}\thanks{
Partially supported by  
KBN grant \mbox{2 P301 020 07}.
On~leave from:
Department of Mathematical Methods in Physics,
Warsaw University, ul.~Ho\.{z}a 74, Warsaw, \mbox{00--682~Poland}.
{\sc e-mail}: {\sc pmh33@amtp.cam.ac.uk} or {\sc pmh@fuw.edu.pl}
}\\
Department of Applied Mathematics and Theoretical Physics,\\
University of Cambridge, Silver St., Cambridge CB3 9EW, England.\\
}
}
\date{DAMTP-97-78}

\maketitle

\vspace*{-7.5mm}
\begin{abstract}
The relationship between the exactness  of a first order differential
calculus on a comodule algebra $P$ and the Galois property of $P$ is investigated.
\end{abstract}

\section*{}\vspace*{-7.5mm}

\hge s describe \qpb s the same way Hopf algebras describe \qg s.
A \hop\ $H$ plays the role of the structure group, and the total space of 
a bundle is replaced by a right $H$-comodule algebra $P$. More precisely,
we have:

\bde\Label{hgdef}
Let $H$ be a Hopf algebra, $P$ a right $H$-comodule algebra with
the coaction $\dr:P\ra P\ot H$, and
$B:=P\sp{co H}:=\{p\in P\,|\; \hD_R\ p=p\ot 1\}$ the subalgebra
of all right coinvariants. We say that
P is a {\em Hopf-Galois extension} (or $H$-Galois extension)
of $B$ iff the canonical left $P$-module right $H$-comodule map
\[
can:=(m\ot id)\circ (id\ot\sb B \hD_R)\, :\; P\ot\sb B P\lra P\ot H,
\]
where $m$ is the multiplication map, is bijective.
\ede
A  \hge\ can also be defined with the help of the universal
differential calculus
(see~\cite{b-t96}, Proposition~1.6 in~\cite{h-p96}):
\bpr\Label{dif}
Let $P$ be a right $H$-comodule algebra and $B$ the subalgebra of
all coinvariants of $P$, as in Definition~\ref{hgdef}. 
Let $\hO^1\! P$, $\hO^1\! B$ and
$H^+:=\mbox{\em Ker}\,\he$ denote
the universal calculus (kernel of the multiplication map) on $P$, $B$
and  the augmentation ideal (kernel of the counit map \he) of $H$, 
respectively. Then  $P$ is
an $H$\gal\ {\em if and only if} the following sequence of left $P$-modules
is exact:
\beq\Label{e}
0\lra P(\hO^1\! B)P\lra\hO^1\! P\st{\overline{can}}{\lra} P\ot H^+\lra 0\, .
\eeq
Here $\overline{can}$ is the appropriate restriction of the
map $(m\ot id)\ci(id\ot\hD_R):P\ot P\ra P\ot H$.
\epr

Recently, it has been observed \cite[Proposition~2.3]{bh97} that the analogous
differential formulation holds true for a more general notion of a
coalgebra-Galois extension. Another direction for broadening the context
of (\ref{e}) is to consider a general differential calculus on~$P$. Let
$N_P\inc\hO^1\! P$ be a $P$-submodule defining a first order right 
covariant ($\dr(N_P)\inc N_P\ot H$) calculus $\hO^1(P):=\hO^1\! P/N_P\,$,
and $R_H\inc H$ an $ad_R$-invariant right ideal defining a bicovariant
differential calculus on $H$~\cite{w-s89}. Then one might demand that
there exists an exact sequence of left $P$-modules:
\beq\Label{e2}
0\lra P\hO^1(B)P\lra\hO^1(P)\st{\chi}{\lra} P\ot (H^+/R_H)\lra 0\, ,
\eeq
where $\chi$ is the quotient map obtained from~$\overline{can}$. 
Sequence (\ref{e2}) is a starting point of the quantum-group gauge theory 
proposed in~\cite{bm93} and continued in~\cite{h-p96}. 

It is natural to expect that the Galois property of an extension
$B\inc P$ (bijectivity of $can$) together with some compatibility 
conditions on $N_P$ and $R_H$ will ensure the existence of exact 
sequence~(\ref{e2}). One can also ask when the existence and exactness
of (\ref{e2}) entails the Galois property for $B\inc P$. The aim of this
note is to provide an answer to this question:

\bpr\Label{main}
Let $H$ be a Hopf algebra and $P$ a right $H$-comodule algebra.
Let $N_P$ be a $P$-subbimodule of
$\hO^1\! P:=\mbox{\em Ker}\, m$, $\op:=\hO^1\! P/N_P$ a first order
differential calculus on $P$, and $V$ a vector subspace of 
$H^+:=\mbox{\em Ker}\,\he$ satisfying the compatibility condition
$\overline{can}(N_P)=P\ot V$. 
Then there exists a map 
\[
\chi:\op\ra P\ot(H^+/V),\;\;\;
\chi([\ha ]_{N_P}):=(id\ot\pi_H)(\chi(\ha)),\;\;\;
 H^+\st{\pi_H}{\ra}H^+/V\, ,
\]
and the following statements are {\em equivalent}:
\be
\item $P$ is an $H$\gal\ of $B:=P^{co H}$. 
\item The sequence of left $P$-modules
\beq\Label{e3}
0\ra P\hO^1(B)P\hookrightarrow\hO^1(P)\st{\chi}{\ra} P\ot (H^+/V)\ra 0
\eeq
is exact and $(N_P\cap\mbox{\em Ker}\,\overline{can})\inc P(\hO^1\! B)P$.
\ee
\epr
\bpf
First observe that we have a commutative diagram (of left $P$-modules) with
exact rows and columns:

\beq\Label{d2}
\def\normalbaselines{\baselineskip24pt
\lineskip3pt \lineskiplimit3pt }
\def\mapright#1{\smash{
\mathop{\!\!\!{-\!\!\!}-\!\!\!\longrightarrow\!\!\!}\limits^{#1}}}
\def\mapdown#1{\Big\downarrow
\rlap{$\vcenter{\hbox{$\scriptstyle#1$}}$}}
\matrix{
0&\mapright{}&\mbox{Ker}\,\overline{\chi}&\mapright{}&\mbox{Ker}\,\overline{can}
&\mapright{}&\mbox{Ker}\,\chi\cr
&&\mapdown{}&&\mapdown{}&&\mapdown{}\cr
0&\mapright{}&N\sb P&\mapright{}
&\Omega\sp 1\! P&\mapright{\pi\sb P}&\op&\mapright{}&0\cr
&&\mapdown{\overline{\chi}}&&\mapdown{\overline{can}}&&\mapdown{\chi}\cr
0&\mapright{}&P\ot V&\mapright{}&P\ot H^+&\mapright{id\ot\pi\sb H}
&P\ot(H^+/V)&\mapright{}&0\cr
&&\mapdown{}&&\mapdown{}&&\mapdown{}\cr
&&0&\mapright{}&\mbox{Coker}\,\overline{can}
&\mapright{}&\mbox{Coker}\,\chi&\mapright{}&0,\cr
}
\eeq 
where $\overline{\chi}$ is the appropriate restriction
of $\overline{can}$. By the Snake Lemma (e.g., see Section~1.2 in~\cite{b-n80}),
we have an exact sequence:
\beq\Label{e4}
0\lra\mbox{Ker}\,\overline{\chi}\lra\mbox{Ker}\,\overline{can}\lra
\mbox{Ker}\,\chi 
\lra 0\lra\mbox{Coker}\,\overline{can}\lra\mbox{Coker}\,\chi\lra 0.
\eeq
Note also that $\mbox{Ker}\,\overline{\chi}=N_P\cap P(\hO^1\! B)P$,
as a consequence of diagram~(\ref{d2}).

Assume now that $(N_P\cap\mbox{Ker}\,\overline{can})\inc P(\hO^1\! B)P$ 
and (\ref{e3})
is exact. Then $\mbox{Coker}\,\chi=0$, and by the exactness of
(\ref{e4}), $\mbox{Coker}\,\overline{can}=\mbox{Coker}\,\chi=0$.
Next, 
as $(N_P\cap\mbox{Ker}\,\overline{can})\inc P(\hO^1\! B)P$
(by assumption) and $P(\hO^1\! B)P\inc\mbox{Ker}\,\overline{can}$
(because $\overline{can}(p(1\ot b-b\ot 1)p')=0$),
we have 
\[
\mbox{Ker}\,\overline{\chi}=N_P\cap\mbox{Ker}\,\overline{can}
=N_P\cap\mbox{Ker}\,\overline{can}\cap P(\hO^1\! B)P=N_P\cap P(\hO^1\! B)P.
\]
On the other hand, we have an exact sequence
\beq\Label{e5}
0 \ra N_P\cap P(\hO^1\! B)P\hookrightarrow P(\hO^1\! B)P \st{\pi}{\ra}
 P\hO^1\!(B)P \ra 0, 
\eeq
where $\pi$ is the canonical surjection. Combining this sequence with the
first part of (\ref{e4}) and taking advantage of $\mbox{Ker}\,\overline{\chi}=
N_P\cap P(\hO^1\! B)P$ and the exactness of (\ref{e3}), we get the following 
commutative diagram with exact rows:

\beq\Label{5}
\def\normalbaselines{\baselineskip30pt
\lineskip3pt \lineskiplimit3pt }
\def\mapright#1{\smash{
\mathop{\lra}
\limits^{#1}}}
\def\mapdown#1{\Big\downarrow
\rlap{$\vcenter{\hbox{$\scriptstyle#1$}}$}}\matrix{
0 &\mapright{}& N_P\cap P(\hO^1\! B)P
&\mapright{\subseteq}& P(\hO^1\! B)P &\mapright{}&
 P\hO^1\!(B)P &\mapright{}& 0 \cr
\mapdown{id}&&\mapdown{id}&&\mapdown{\subseteq}&&\mapdown{id}&&\mapdown{id}\cr
0 &\mapright{}& N_P\cap P(\hO^1\! B)P
&\mapright{\subseteq}& \mbox{Ker}\,\overline{can}
&\mapright{}& P\hO^1\!(B)P &\mapright{}& \phantom{id.}0\phantom{id}. \cr
}
\eeq\ \\ \ \\
Therefore, by the Five Isomorphisms Lemma (e.g., 
see~\cite[Section~1.2, Corollary~3]{b-n80}), 
$P(\hO^1\! B)P=\mbox{Ker}\,\overline{can}$. Thus we have shown that
$\mbox{Coker}\,\overline{can}=0$ and $\mbox{Ker}\,\overline{can}=P(\hO^1\! B)P$.
Now it follows from Proposition~\ref{dif} that $P$ is an $H$\gal.

Conversely, assume that $P$ is an $H$\gal. Then again by Proposition~\ref{dif},
$\mbox{Ker}\,\overline{can}=P(\hO^1\! B)P$ and $\mbox{Coker}\,\overline{can}=0$.
With the help of (\ref{e4}), the latter implies that 
$\mbox{Coker}\,\chi=0$. As $\mbox{Ker}\,\overline{can}=P(\hO^1\! B)P$,
the first part of (\ref{e4}) can be written as the following exact sequence:
\beq\Label{e6}
0 \ra N_P\cap P(\hO^1\! B)P\hookrightarrow P(\hO^1\! B)P \st{\pi}{\ra}
\mbox{Ker}\,\chi \ra 0, 
\eeq
so that $\mbox{Ker}\,\chi=P\hO^1(B)P$. Thus we have shown that
(\ref{e3}) is exact. To complete the proof it suffices to note that
$N_P\cap\mbox{Ker}\,\overline{can}=N_P\cap P(\hO^1\! B)P\subseteq P(\hO^1\! B)P$. 
\epf\ \\

Note that to give (\ref{e3}) a more differential-geometric meaning, one
assumes that $V=R_H$ is an $ad_R$-invariant right ideal of $H$, so that
the corresponding differential calculus on $H$ is bicovariant. For the
same reason, one also assumes that 
 $\dr(N_P)\inc N_P\ot H$, so that $\op:=\hO^1\! P/N_P$
is right covariant. These particulars are not assumed in Proposition~\ref{main}
as they are unnecessary to prove it. To end this note,
let us provide an alternative to Proposition~\ref{main}:

\bpr\Label{main2}
Let $H$ be a Hopf algebra and $P$ a right $H$-comodule algebra 
such that the canonical map $P\ot_BP\st{can}{\lra}P\ot H$ is surjective. Let 
$\op:=\hO^1\! P/N_P$ be a first order
differential calculus on~$P$, and $V$ a vector subspace of 
$H^+:=\mbox{\em Ker}\,\he$ satisfying the compatibility condition
$\overline{can}(N_P)\inc P\ot V$. 
Then the following statements are {\em equivalent}:
\be
\item $P$ is an $H$\gal\ of $B:=P^{co H}$ and 
      $P\ot V\inc\overline{can}(N_P)$. 
\item The sequence of left $P$-modules
\beq\Label{e9}
0\ra P\hO^1(B)P\hookrightarrow\hO^1(P)\st{\chi}{\ra} P\ot (H^+/V)\ra 0,
\eeq
where $\chi$ is defined as in Proposition~\ref{main},
is exact and $(N_P\cap\mbox{\em Ker}\,\overline{can})\inc P(\hO^1\! B)P$.
\ee
\epr
\bpf
This proof is very similar to the proof of Proposition~\ref{main}. Note first
that the surjectivity of $can$ implies the surjectivity of
$(m\ot id)\ci(id\ot\hD_R):P\ot P\ra P\ot H$.
With the help of~\cite[(3)]{h-p96}, this in turn implies 
$\mbox{Coker}\,\overline{can}=0$.
Again, we can take advantage of the Snake Lemma to obtain an exact sequence
\beq\Label{e7}
0\lra\mbox{Ker}\,\overline{\chi}\lra\mbox{Ker}\,\overline{can}\lra
\mbox{Ker}\,\chi 
\lra\mbox{Coker}\,\overline{\chi}\lra 0\lra\mbox{Coker}\,\chi\lra 0.
\eeq

Assume first that (\ref{e9}) is exact and 
$(N_P\cap\mbox{Ker}\,\overline{can})\inc P(\hO^1\! B)P$. It follows from the 
exactness of (\ref{e9}) that $\mbox{Ker}\,\chi=P\ob P$. Hence,
as $P(\hO^1\! B)P\inc\mbox{Ker}\,\overline{can}$, the canonical quotient
map $\mbox{Ker}\,\overline{can}\ra\mbox{Ker}\,\chi=P\ob P$ of
sequence (\ref{e7}) is necessarily surjective. Consequently, from the
exactness of~(\ref{e7}), we can infer that $\mbox{Coker}\,\overline{\chi}=0$,
i.e., $\overline{can}(N_P)=P\ot V$. 
We already know that $\mbox{Coker}\,\overline{can}=0$
and can now proceed as in the proof of Proposition~\ref{main} to show that
$P$ is an $H$\gal.

Conversely, assume that $P$ is an $H$\gal\ of $B:=P^{co H}$ and 
$P\ot V\inc\overline{can}(N_P)$. The latter implies that 
$\mbox{Coker}\,\overline{\chi}=0$.
Now one can argue as in the proof of Proposition~\ref{main} to conclude that
(\ref{e9}) is exact and $(N_P\cap\mbox{Ker}\,\overline{can})\inc P(\hO^1\! B)P$.
\epf

\section*{Acknowledgements}
It is a pleasure to thank Tomasz Brzezi\'nski for a helpful discussion.

\end{document}